\input amstex
\documentstyle{amsppt}
\input epsf.tex
\pageheight{47pc}

\def\nologo{\let\logo@\relax}
\def\stmin{\smallsetminus}

\def\d{\partial}
\def\conj{\text{\rm conj}}
\def\ge{\geq}
\def\le{\leq}
\def\Cl{\roman{Cl}}
\def\Int{\roman{Int}}
\def\D{\Delta}
\def\G{\Gamma}
\def\b{\beta}
\def\la{\langle}
\def\ra{\rangle}
\def\trans{\pitchfork}
\def\vbot{\bot\!\!\!\bot}
\def\R{\Bbb{R}}
\def\Z{\Bbb{Z}}
\def\C{\Bbb{C}}
\def\Cp#1{\Bbb{C}\roman{P}^{#1}}
\def\Rp#1{\Bbb{R}\roman{P}^{#1}}
\nologo
\magnification=1200
\NoRunningHeads
\NoBlackBoxes

\topmatter
\title
Knotting of algebraic curves in $\Cp2$
\endtitle
\author Sergey Finashin \endauthor
\address
Middle East Technical University,
Ankara 06531 Turkey
\endaddress
\email
serge\,\@\,metu.edu.tr
\endemail
\abstract
For any $k\ge 3$, I construct infinitely many pairwise smoothly non-isotopic
smooth surfaces $F\subset\Cp2$ homeomorphic to a non-singular
algebraic curve of degree $2k$, realizing the same homology class
as such a curve
and having  abelian fundamental group $\pi_1(\Cp2\stmin F)$.
  This gives an answer to Problem 4.110 in the Kirby list \cite{K}.
\endabstract
\endtopmatter

\document

\heading
1. Introduction 
\endheading
The goal of this paper is the following result.

\proclaim{1.1. Theorem}
For any $k\ge3$ there exist infinitely many 
smooth oriented closed surfaces
$F\subset\Cp2$ representing class $2k\in H_2(\Cp2)=\Z$,
having $\roman{genus}(F)=(k-1)(2k-1)$
and $\pi_1(\Cp2\stmin F)\cong\Z/2k$, such that
the pairs $(\Cp2,F)$ are pairwise smoothly non-equivalent.
\endproclaim

Here smooth non-equivalence of $(\Cp2,F)$ and $(\Cp2,F')$ means that
$F$ cannot be transformed into $F'$ by a diffeomorphism
$\Cp2\to\Cp2$.

Construction of these surfaces is based on a modification of the
"rim-surgery" of Fintushel and Stern \cite{FS2}, which is applied
for knotting a surface  along an annulus membrane.
By an annulus membrane for
a smooth surface $F$ in a $4$-manifold $X$ I mean a smoothly
embedded surface $M\subset X$, $M\cong S^1\times I$, with
$M\cap F=\d M$ and such that $M$ comes to $F$ normally along $\d M$.
 Assume that such a membrane has
framing $0$, or equivalently, admits a diffeomorphism
of its regular neighborhood $\phi\:U\to S^1\times D^3$
 mapping $U\cap F$ onto $S^1\times f$,
where $f=I\vbot I\subset D^3$ is a disjoint union of two segments,
which are unknotted and unlinked in $D^3$,
that is to say that
 a union of $f$ with a pair of arcs on a sphere $\d D^3$
bounds a trivially
embedded band, $b\subset D^3$, $b\cong I\times I$,  so that
$f= I\times(\d I)\subset b$ (see Figure~1).
The annulus $M$ can be viewed
as $S^1\times\{\frac12\}\times I$ in
$S^1\times b\subset S^1\times D^3\cong U$.

If $X$ and $F$ are oriented, then $f$ inherits an orientation
as a transverse intersection, $f=F\trans D^3$, and
we may choose a band $b$ so that the orientation of $f$
is induced from some orientation of $b$.
It is convenient to view $f=I\vbot I$ as is shown on Figure 1,
so that the segments of $f$ are parallel and oppositely oriented,
with $b$ being a thin band between them.
Such a presentation is always possible if we allow a modification
of $\phi$, since one of the segments of $f$ may be turned around
by a diffeomorphism of $D^3\to D^3$ preserving the other segment fixed.

 Given a knot $K\subset S^3$, we construct a new smooth surface,
 $F_{K,\phi}$, obtained from $F$ by tying a pair of segments
$I\vbot I$ along $K$ inside $D^3$, as is shown on Figure~1.
More precisely, we consider a band $b_K\subset D^3$ obtained from
$b$ by knotting along $K$
and let $f_K$ denote the pair of arcs bounding $b_K$ inside $D^3$.
We assume that the framing of
$b_K$ is chosen the same as the framing of $b$, or equivalently,
that the inclusion homomorphisms from
$H_1(\d D^3\stmin (\d f))=H_1(\d D^3\stmin (\d f_K))$
to $H_1(B^3\stmin f)$ and to $H_1(B^3\stmin f_K)$ have the same kernel.
Then $F_{K,\phi}$ is obtained from $F$ by replacing
$S^1\times f\subset S^1\times D^3\cong U$ with
$S^1\times f_K$.
It is obvious that
$F_{K,\phi}$ is homeomorphic to $F$ and realizes
the same homology class in $H_2(X)$.

\midinsert
\epsfbox{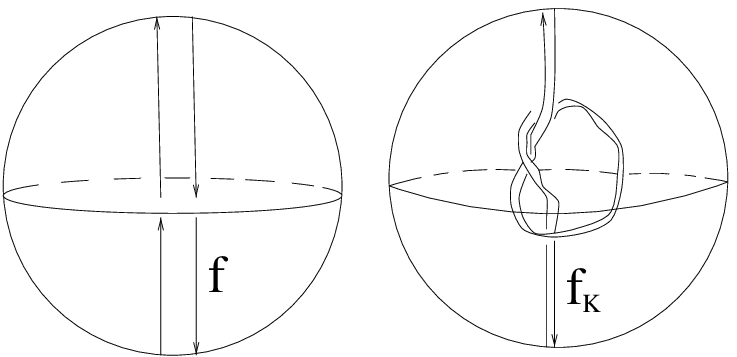}
\botcaption{Figure 1}
Knotting of a band $b_K$
\endcaption
\endinsert

The above construction is called in what follows
{\it an annulus rim-surgery}, since it looks like the
 rim-surgery of Fintushel and Stern \cite{FS2}, except that
we tie two knits together instead of one.
Recall that the usual rim-surgery is applied in \cite{FS2} to
surfaces $F\subset X$ which are {\it primitively embedded},
 that is $\pi_1(X\stmin F)=0$, which is not the case for
the algebraic curves in $\Cp2$ of degree $>1$.
The primitivity condition  is required to
preserve the fundamental group of $X\stmin F$ throughout the knotting.
 An annulus rim-surgery may preserve {\it a non-trivial}
group $\pi_1(X\stmin F)$, if we require
commutativity of $\pi_1(X\stmin(F\cup M))$, instead of primitivity
of the embedding.

\proclaim{1.2. Proposition}
Assume that $X$ is a simply connected closed $4$-manifold,
$F\subset X$ is an oriented closed surface with an annulus-membrane
$M$ of index $0$, $\phi\: U\to S^1\times D^3$ is a trivialization like
described above and $K\subset S^3$ is any knot.
Assume furthermore that
 $F\stmin\d M$ is connected and the group $\pi_1(X\stmin(F\cup M))$
is abelian.
Then the group $\pi_1(X\stmin F_{K,\phi})$ is cyclic and isomorphic to
$\pi_1(X\stmin F)$.
\endproclaim

To prove Theorem 1.1, I apply an annulus rim-surgery
for $X=\Cp2$ letting $F=\C A$ be the complex point set of
a suitable non-singular real algebraic curve,
containing an annulus, $M$, among
the connected components of $\Rp2\stmin\R A$,
where $\R A=\C A\cap\Rp2$ is the real locus of the curve.

One may take, for instance, a real algebraic curve
$\C A$ of degree $d=2k$, with $\R A$
containing $k$ components (called ovals), $O_1,\dots,O_k$,
such that $O_i$ lies inside $O_{i+1}$ in $\Rp2$, $i=1,\dots,k-1$.
 Such a real algebraic curve, known as a {\it maximal nest} curve,
can be constructed by a small perturbation of
a union of $k$ real conics, whose real parts (ellipses) are
ordered by inclusion in $\Rp2$.
The connected components of $\Rp2\stmin\R A$ include a disc, $R_0$,
inside $O_1$, a M\"obius band component, $R_k$, outside $O_k$,
and $k-1$ annuli, $R_i$, between $O_i$ and $O_{i+1}$, whose closures,
$\Cl(R_i)$ are obviously $0$-framed annulus-membranes on $\C A$.
For simplicity, let us choose $M=\Cl(R_1)$.

\proclaim{1.3. Proposition}
The assumptions of Proposition 1.2 hold for $X=\Cp2$,
$F=\C A$ being a maximal nest real algebraic curve of degree
 $2k\ge6$ and $M=\Cl(R_1)$.
\endproclaim

Assuming that the class $[F]\in H_2(X;\Z/2)$ vanishes,
one can consider a double covering $p\:Y\to X$ branched along $F$;
such a covering is unique if we require in addition that $H_1(X;\Z/2)=0$.
 Similarly, we consider the double coverings
$Y(K,\phi)\to X$ branched along $F_{K,\phi}$.
To prove non-equivalence of pairs $(\Cp2,F_{K,\phi})$
for some family of knots $K$,
it is enough to show that $Y(K,\phi)$ are not pairwise diffeomorphic.
To show it,
I use that $Y(K,\phi)$ is diffeomorphic to the $4$-manifolds
$Y_{K\#K}$ obtained from $Y$
by a surgery introduced in \cite{FS1}
({\it FS-surgery}).

\proclaim{1.4. Proposition}
$Y(K,\phi)$ is diffeomorphic to
a  $4$-manifold obtained from $Y$ by the FS-surgery along
the torus $T=p^{-1}(M)$ via the knot $K\# K\subset S^3$.
\endproclaim

To distinguish the diffeomorphism types of $Y_{K\#K}$
one can use the
 formula of Fintushel and Stern \cite{FS1} for SW-invariants of
a $4$-manifold $Y$ after FS-surgery along a torus $T\subset Y$.
Recall that this formula can be applied if the
SW-invariants of $Y$ are well-defined and a torus $T$,
realizing a non-trivial class $[T]\in H_2(Y)$, is {\it c-embedded}
(the latter means that $T$ lies as a non-singular fiber
in a cusp-neighborhood in $Y$, cf. \cite{FS1}).
Being an algebraic surface of genus $\ge1$, the double plane
$Y$ has well-defined
SW-invariants. The conditions on $T$ are also satisfied.

\proclaim{1.5. Proposition}
Assume that $X$, $F$ and $M$ are like in Proposition 1.2.
Then the torus $T=p^{-1}(M)$ is primitively embedded in $Y$
and therefore $[T]\in H_2(Y)$ is an infinite order class.
If, moreover, $X$, $F$ and $M$ are like in Proposition 1.3, then
$T\subset Y$ is c-embedded.
\endproclaim

Recall that the product formula \cite{FS1}
$$SW_{Y_K}=SW_Y{}^{.}\D_K(t),\ \text{ where } t=\exp(2[T])
$$
expresses the Seiberg-Witten invariants
(combined in a single polynomial) of the manifold $Y_K$, obtained
by an FS-surgery, in terms of the Seiberg-Witten invariants of $Y$
and the Alexander polynomial, $\D_K(t)$, of $K$.

 This formula implies that the basic classes of $Y_K$ can be expressed
as $\pm \b+ 2n[T]$, where
$\pm\b\in H_2(Y)$ are the basic classes of $Y$ and
 $|n|\le\deg(\D_K(t))$,
are the degrees of the non-vanishing monomials in $\D_K(t)$.
By the adjunction formula, the class $[T]$ is orthogonal
to all the basic classes, $\b$, of $Y$.
So, if $[T]$ has infinite order,
then the manifolds $Y(K,\phi)\cong Y_{K\#K}$ differ from each other by 
their SW-invariants,
and moreover, just by the numbers of their basic classes,
for a suitable infinite family of knots $K$,
since the number of the basic classes is determined by
the number of the terms in $\D_{K\# K}=(\D_K)^2$
(one can take for instance a family $\{K_i\}$, where $K_i$
is a connected sum of $i$ copies of a trefoil knot).

\head 2. Commutativity of the fundamental group throughout
the knotting
\endhead

\proclaim{2.1. Lemma}
The assumptions of Proposition 2.1 imply that
$\pi_1(X\stmin(F\cup M))=\pi_1(X\stmin F)$
is cyclic with a generator presented
by a loop around $F$.
\endproclaim

\demo{Proof}
The Alexander duality in $X$ combined with the exact cohomology
sequence of a pair $(X,F\cup M)$ gives
$$
H_1(X\stmin(F\cup M))\cong H^3(X,F\cup M)=H^2(F\cup M)/i^*
H^2(X)
$$
where $i\: F\cup M\to X$ is the inclusion map.
 If $F$ is oriented and $F\stmin\d M$ is connected, then
the Mayer-Vietoris Theorem yields
$H^2(F\cup M)\cong H^2(F)\cong\Z$, and thus
$H_1(X\stmin(F\cup M))\cong H_1(X\stmin F)$ is cyclic
with a generator presented by a loop around $F$.
The same property holds for  the fundamental groups of
$X\stmin(F\cup M)$ and $X\stmin F$,
since they are abelian by the assumption of Proposition 2.1.
\qed
\enddemo

\demo{Proof of Proposition 1.2}
Put $X_0=\Cl(X\stmin U)$. Then $\d X_0=\d U\cong S^1\times S^2$ and
$\d U\stmin F$ is a deformational retract of
$U\stmin(F\cup M)$, so
$$
\pi_1(X_0\stmin F)=\pi_1(X\stmin(F\cup M))
$$
Since this group is cyclic and is generated by a loop around $F$,
the inclusion homomorphism
$h\:\pi_1(\d U\stmin F)\to\pi_1(X_0\stmin F)$
is epimorphic and thus
$\pi_1(X_0\stmin F)=\pi_1(\d U\stmin F)/k$, where $k$
is the kernel of $h$.

Applying the Van Kampen theorem to the triad
$(X_0\stmin F,U\stmin F_{K,\phi},\d U\stmin F)$,
we conclude that
$
\pi_1(X\stmin F_{K,\phi})\cong\pi_1(U\stmin F_{K,\phi})/j(k)
$,
where $j\:\pi_1(\d U\stmin F)\to\pi_1(U\stmin  F_{K,\phi})$
is the inclusion homomorphism.
Furthermore, in the splitting
$$
\pi_1(U\stmin  F_{K,\phi})\cong\pi_1(S^1\times(D^3\stmin
f_K))\cong\Z\times\pi_1(D^3\stmin f_K)
$$
factorization by $j(k)$ kills the first factor $\Z$ and adds some
relations to $\pi_1(D^3\stmin f_K)$, one of which
effects to $\pi_1(D^3\stmin f_K)$ as if we attach
a $2$-cell  along a loop, $m_b$, turning
around the band $b_K$
(to see it, note that factorization by $k$ leaves only one generator of
$\pi_1(\d D^3\stmin f_K)=\pi_1(S^2\stmin\{4 \text{\rm{pts}}\})$).
Attaching such a $2$-cell effects to $\pi_1$
as connecting together a pair of the endpoints
of $f_K$, which transforms $f_K$ into an arc (see Figure~2).
This arc is unknotted and thus factorization by $j(k)$ makes
$\pi_1(D^3\stmin f_K)$ cyclic and leaves
$\pi_1(X\stmin F_{K,\phi})$ isomorphic to
$
\pi_1(X_0\stmin F)\cong\pi_1(X\stmin(F\cup
M))\cong\pi_1(X\stmin F)
$.
\qed
\enddemo

\midinsert
\epsfbox{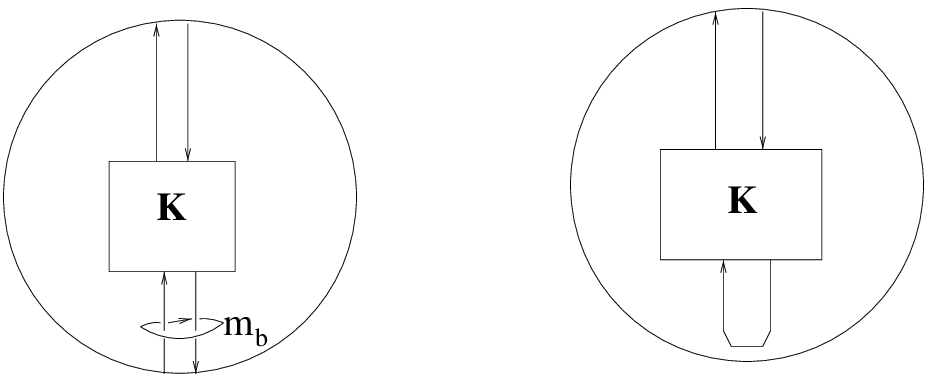}
\botcaption{Figure 2}
Gluing a $2$-cell along $m_b$ effects as transforming
 $f_K$ into an unknotted arc
\endcaption
\endinsert

\demo{Proof of Proposition 1.3}
All the assumptions of Proposition 1.2 except the last two are
obviously satisfied. It is well known
that $\C A\stmin\R A$ splits for a maximal nest curve $\C A$
into a pair of connected components permuted by the complex conjugation,
and thus, $\C A\stmin\d M$ is connected, provided $\d M\subsetneqq\R A$,
which is the case for $k\ge 3$
(recall in turn that $\C A\stmin\R A$ has not more then $2$
connected components
for any non-singular real algebraic curve $\C A\subset\Cp2$).
So, it is left to check only the last assumption,
that the group $\pi_1(\Cp2\stmin(\C A\cup M))$ is abelian.
This follows from an explicit description of the homotopy type
and, in particular, the fundamental group of the complement
$\Cp2\stmin(\C A\cup\Rp2)$ given in \cite{F1} (see also \cite{FKV,\S4})
for $L$-curves.
Recall that a non-singular curve $\C A\subset\Cp2$ of degree $m$
is an $L$-curve if it can be obtained by a non-singular
perturbation from a curve $\C A_0=\C L_1\cup\dots\C L_m$
splitting into $m$ real lines, $\C L_i$, in a generic position.
The maximal nest curves, $\C A\subset\Cp2$, are known \cite{F1} to
be $L$-curves with the presentation
$\pi=\pi_1(\Cp2\stmin(\C A\cup\Rp2))=\la a,b\,|\,a^{2k}b^{2k}=1\ra$,
where $a$, $b$ are represented by loops around
the two connected components of $\C A\stmin\R A$;
a basis point and these loops are taken on the conic
$C=\{x^2+y^2+z^2=0\}\subset\Cp2$, which have the real point set empty,
cf. \cite{F1}.
%
The group $\pi_1(\Cp2\stmin(\C A\cup M))$ is obtained from $\pi$
by adding the relations corresponding to
puncturing the components $R_i$, $0\le i\le k$, $i\ne 1$,
of $\Rp2\stmin\R A$.
Such a relation (for $R_i$) is $a^{2k-i}b^i=b^{2k-i}a^i=1$,
which follows from the results of \cite{FKV,\S4}
(see also \cite{F2} for detailed proofs).
A pair of the relations for $i=2$ and $i=3$ implies that $a=b$.

For convenience of a reader, I included in the appendix
a brief review of the arguments in \cite{F2} and \cite{FKV,\S4}
relevant to the above calculation.
\qed
\enddemo

\remark{Remarks}
\roster
\item
It is well-known that a real curve with a maximal nest arrangement of ovals is
unique up to a rigid isotopy, that is an isotopy in the class of
non-singular real curves.
\item
It follows from the proof above that $\pi_1(\Cp2\stmin(\C A\cup M))$
is not abelian and $\C A\stmin\d M$ is not connected
for a maximal nest quartic, $\C A$.
\item
The other types of real algebraic curves
containing a pair of ovals, $O_1$, $O_2$,
bounding an annulus, may be equally good for the knotting
construction, at least if they are $L$-curves.
The computation of $\pi$ may become even more simple.
For instance, for any $L$-curve  $\C A$
of degree $2k\ge6$ whose real part $\R A$
contains {\it just two ovals}, $O_1$, $O_2$
(no matter if they bound an annulus or lie separately),
we have
$\pi_1(\Cp2\stmin(\C A\cup\Rp2))\cong\Z/4k$.
\endroster
\endremark

\head 3. The double surgery in the double covering
\endhead

\demo{Proof of Proposition 1.4}
 The proof is based on the following two observations.
First, we notice that $Y(K,\phi)$ is obtained from $Y$
by a pair of FS-surgeries along the tori parallel to $T$,
then we notice that such pair of surgeries is equivalent
to a single FS-surgery along $T$.
The both observations are corollaries of \cite{FS2, Lemma 2.1},
so, I have to recall first the construction from \cite{FS1, FS2}.

An FS-surgery \cite{FS1} on a $4$-manifold $X$ along a torus
$T\subset X$, with the self-intersection
$T\circ T=0$, via a knot $K\subset S^3$ is defined as
a {\it fiber sum} $X\#_{T= S^1\times m_K} S^1\times M_K$,
that is an amalgamated connected sum of $X$ and $S^1\times M_K$
along the tori $T$ and $S^1\times m_K\subset  S^1\times M_K$.
Here $M_K$ is a $3$-manifold obtained by the $0$-surgery along $K$ in $S^3$,
and $m_K$ denotes a meridian of $K$ (which may be seen both in
$S^3$ and in $M_K$).
 Such a fiber sum operation can be viewed as a direct product
of $S^1$ and the corresponding $3$-dimensional operation,
which I call {\it $S^1$-fiber sum}.

More precisely, $S^1$-fiber sum
$X\#_{K=L}Y$ of oriented $3$-manifolds
$X$ and $Y$ along oriented framed knots $K\subset X$ and $L\subset Y$
 is the manifold obtained by gluing the complements
$\Cl(X\stmin N(K))$ and $\Cl(Y\stmin N(L))$ of tubular neighborhoods,
$N(K)$, $N(L)$, of $K$ and $L$ via a diffeomorphism
$f\:\d N(K)\to\d N(L)$ which identifies the longitudes
of $K$ with the longitudes of $L$ preserving their orientations,
and the meridians of $K$ with the meridians of $L$ reversing the
orientations.
 As it is shown in \cite{FS2, Lemma 2.1}, tying
an arc in $D^3$ along a knot $K\subset S^3$
can be interpreted as a fiber sum
$D^3\#_{m=m_K}M_K$, where $m$ is a meridian around this arc.
The meridians $m$ and $m_K$
are endowed here with the $0$-framings ($0$-framing of a meridian
makes sense as a meridian lies in a small $3$-disc).
 To understand this observation, it is useful to view
an $S^1$-fiber sum with $M_K$
as surgering a tubular neighborhood, $N(m)$, of $m$
and replacing it by the complement,
 $S^3\stmin N(K)$ of a tubular
neighborhood, $N(K)$, of $K$, so that the longitudes of $m$ are
glued to the meridians of $K$ and the meridians of $m$ to the
longitudes of $K$.
The framing of an arc in $D^3$ is preserved under such a fiber  sum,
so tying the band $b\subset D^3$ is equivalent to taking
an $S^1$-fiber sum with $M_K$ along
a meridian $m_b$ around $b$.

The double covering over $D^3$ branched along $f$
is a solid torus, $N\cong S^1\times D^2$, and the pull back of $m_b$
splits into a pair of circles, $m_1,m_2\subset N$, parallel to
$m=S^1\times\{0\}$.
Therefore,
$Y(K,\phi)$ is obtained from $Y$ by performing FS-surgery twice,
along the tori
$$
T_i=S^1\times m_i\subset p^{-1}(U)\cong S^1\times N,\ \
i=1,2
$$
The following Lemma implies that this gives the same result as
a single FS-surgery along $T=p^{-1}(M)$
via the knot $K\#K$.
\qed
\enddemo

\proclaim{3.1. Lemma}
For any pair of knots, $K_1, K_2$, the manifold
$$M_{K_1}\#_{m_{K_1}=m_1}
N\#_{m_2=m_{K_2}}M_{K_2}$$
obtained by taking an $S^1$-fiber sum twice,
is diffeomorphic to
$N\#_{m=m_{K}}M_{K}$, for  $K=K_1\# K_2$,
via a diffeomorphism identical on $\d N$.
\endproclaim

\demo{Proof}
A solid torus $N$ can be viewed as the complement $N=S^3-N'$ of an
open tubular neighborhood $N'$ of an unknot, so that $m,m_1,m_2$ represent
meridians of this unknot. Taking a fiber sum of $S^3$ with
 $M_{K_i}$ along $m_i=m_{K_i}$
is equivalent to knotting  $N'$ in $S^3$ via $K_i$.
So, performing
 $S^1$-fiber sum twice, along $m_1$ and $m_2$, we obtain the same result
as after taking fiber sum along $m$ once, via $K=K_1\# K_2$.
\qed\enddemo

\remark{Remark}
The above additivity property can be equivalently stated as
$$
M_{K_1}\#_{m_{K_1}=m_{K_2}}M_{K_2}\cong M_{K_1\#K_2}
$$
\endremark

\demo{Proof of Proposition 1.5}
 Lemma 2.1 implies that, in the assumptions of Proposition 1.2,
$\pi_1(Y\stmin(F\cup T))$ is a cyclic group with
 a generator represented by a loop around $F$.
Thus, $\pi_1(Y\stmin T)=0$ and, by the Alexander duality,
$ H_3(Y,T)=H^1(Y\stmin T)=0$, which implies that
$[T]\in H_2(Y)$ has infinite order.

To check that $T$ is c-embedded it is enough to observe that
there exists a pair of vanishing cycles on $T$,
or more precisely, a pair of
$D^2$-membranes, $D_1,D_2\subset Y$, on $T$,
having $(-1)$-framing and intersecting at
a unique point $x\in T$, so that
$[\d D_1],[\d D_2]$ form a basis of $ H_1(T)$.
In the setting of Proposition 1.3,
$Y$ is the complex point set of a real
algebraic surface and $T$ is a connected component of its real part.
Nodal degenerations of a curve $\C A$ shown on Figure 3
give nodal degenerations of the double covering, $Y$, which
provide the membranes $D_i$. These $D_i$ are
the halves of the $(-2)$-spheres in $Y$, i.e., of the corresponding
complex vanishing cycles.
One of these complex cycles is just $p^{-1}(R_0)$; it
corresponds to contraction of the oval $O_1$ into a node.
The other complex cycle similarly corresponds to a nodal degeneration
of $\C A$ fusing $O_1$ and $O_2$.
\qed\enddemo

\midinsert
\epsfbox{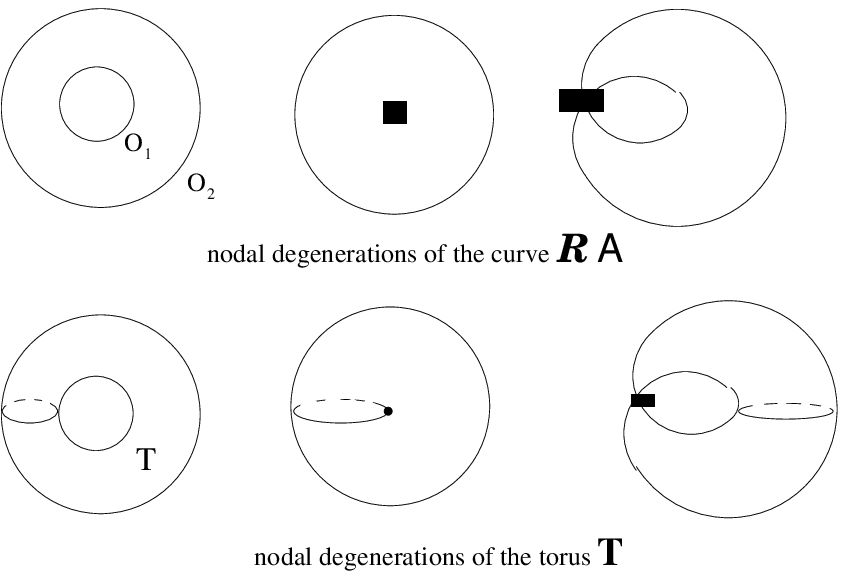}
\botcaption{Figure 3}
Nodal degenerations of $\R A$ providing
$(-1)$-framed $D^2$-membranes on $T$
\endcaption
\endinsert

\comment
\head
4. Concluding remarks
\endhead

The idea behind the construction of the surfaces $F_{K,\phi}$
may resemble the idea \cite{FKV} of producing exotic knottings
of non-orientable surfaces in $S^4$ by factorization of
the Dolgachev surfaces (being
the exotic versions of an elliptic surface $E(1)$)
by the complex conjugation.
 The knotting construction of \cite{FKV} was found as
the quotient of a logarithmic transform equivariant with respect
to the complex conjugation.
In a similar sense,
the annulus rim-surgery appears to be the quotient
 of the Fintushel-Stern surgery \cite{FS1}.
 This analogy may look not complete, because the quotient map
$Y\to\Cp2$ in this paper is holomorphic, whereas in \cite{FKV},
it is used the quotient by an anti-holomorphic involution.
However,  it is interesting to note
that for a K3 surface, $X$, the distinction between holomorphic
and anti-holomorphic vanishes. More precisely,
any anti-holomorphic involution turns out to be holomorphic
with respect to another complex structure on $X$ (cf. \cite{D}).

 Furthermore, the statement of Theorem 1.1 in the special case
 $k=3$ can be proved in another way, using the technique of \cite{FKV}
applying complex-conjugation-equivariant logarithmic transforms
to a K3 surface (viewed as a real elliptic surface $E(2)$)
 in the same way as it was done in \cite{FKV} for $E(1)$.
More generally, taking quotients by the complex conjugation
of various real forms of K3 surfaces, we obtain
rational surfaces with an  anti-bi-canonical curves as the branching
loci. Applying logarithmic transforms to the K3 surfaces, we obtain
exotic knottings of these curves.
\endcomment

\head  Appendix: The topology of
$\Cp2\stmin(\Rp2\cup\C A)$ for $L$-curves $\C A$
\endhead
\eightrm

Let $\C A_0=\C L_1\cup\dots\cup\C L_m\subset\Cp2$ denote 
the complex point set of a real curve of
degree $m$ splitting into $m$ lines, $\C L_i$.
Put $\widetilde V=C\cap\C A_0$, where $C$ is the conic from
the proof of Proposition 1.3.
 Our  first observation is that $C\stmin \widetilde V$ is 
 a deformational retract of
 $\Cp2\stmin(\Rp2\cup\C A_0)$,
and moreover, the latter complement is homeomorphic to
$(C\stmin \widetilde V)\times \Int(D^2)$.
 To see it, it suffices to note that $\Cp2\stmin\Rp2$ is fibered
 over $C$ with a $2$-disc fiber, each fiber being
a real semiline, that is a connected component
of $\C L\stmin\R L$ for some real line $\C L\subset\Cp2$,
where $\R L=\C L\cap\Rp2$.
This fibering maps
a semiline into its intersection point with $C$.

It is convenient to view
the quotient $C/\conj$ of the conic $C$ by the complex conjugation
as the projective plane, $\widehat{\Rp{}}^2$,
dual to $\Rp2\subset\Cp2$, since each real line, $\C L$,
intersects $C$ in a pair of conjugated points.
If we let
 $V=\{l_1,\dots,l_m\}\subset\widehat{\Rp{}}^2 $ denote the set of points
$l_i$ dual to the lines $\R L_i\subset\Rp2$, then
$\widetilde V=q^{-1}(V)$, where 
$q\: C\to C/\conj$ is the quotient map.

The information about a perturbation of $\C A_0$ is encoded in a
{\it genetic graph of a perturbation}, $\G\subset\widehat{\Rp{}}^2$.
The graph $\G$ is a complete graph with the vertex set $V$,
whose edges are  line segments.
Note that there exist two topologically distinct
perturbations of a real node of $\R A_0$ at $p_{ij}=\R L_i\cap\R L_j$,
as well as there exist
two line segments in $\widehat{\Rp{}}^2$ connecting
vertices $l_i, l_j\in V$.
Let $\R A$ denotes a real curve
obtained from $\R A_0$ by a suffucuently small perturbation.
Then the edge of $\G$ connecting $l_i$ and $l_j$ 
contains the points dual to those lines passing through $p_{i,j}$
which do not intersect $\R A$ locally, in a small neighborhood of $p_{i,j}$.

\midinsert
\epsfbox{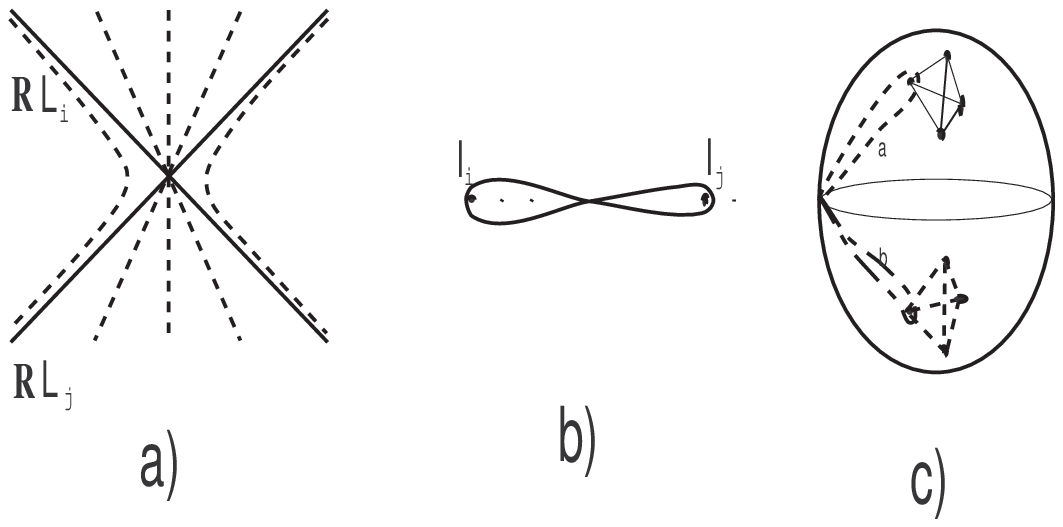}
\botcaption{Figure 4}
a) A perturbation of a real node; the dashed lines are dual
to the points of an edge of $\G$;
b) A figure-eight loop along an edge of $\G$;
c) The loops in $C\stmin \widetilde V$ representing
generators ``a'' and ``b'' 
\endcaption
\endinsert

The complement $\Cp2\stmin(\C A\cup\Rp2)$ turns out to be
homotopy equivalent to a $2$-complex obtained from
$C\stmin \widetilde V$ by adding $2$-cells glued along a figure-eight
shaped loops along the edges of $\widetilde{\G}=q^{-1}(\G)\subset C$.
Such $2$-cells identify pairwise certain
generators of $\pi_1(C\stmin\widetilde V)$
``along the edges'' of $\widetilde{\G}$ (cf. \cite{FKV} for details).
This easily implies that the group
$\pi_1(\Cp2\stmin(\C A\cup\Rp2))$ is generated by a pair of elements,
$a$ and $b$, represented by a pair of loops in $C\stmin \widetilde V$ around
a pair of conjugated vertices of $\widetilde V$.

For example, for a maximal nest curve, the graph $\G$ is contained
in an affine part of $\widehat{\Rp{}}^2$, i.e., has no common points
with some line in $\widehat{\Rp{}}^2$, namely, with a line
dual to a point inside the inner oval of the nest.
Therefore, the graph $\widetilde{\G}$ splits into two connected components
separated by a big circle in $C$.
 A loop around any vertex of $\widetilde V$ from one of these components
 represents $a$, and a loop around a vertex from the other
 component represents $b$.
It is trivial to observe also the relation $a^mb^m=1$
(which is indeed a unique relation in the case of maximal nest curves).

As we puncture $\Rp2$ at a point $x\in \Rp2\stmin\R A_0$,
we attach a $2$-cell to $C\stmin \widetilde V$ along the big circle 
$S_x\subset C$ dual to $x$.
If $x$ moves across a line $\R L_i$, then $S_x$ moves across the pair
of points $q^{-1}(l_i)$.
 Since a small perturbation and puncturing are located at
distinct points of $\Cp2$ and can be done independently,
it is not difficult to see that if we choose $x\in R_i$ 
(in the case of a maximal nest curve $\C A$), then the big circle $S_x$
cuts $C$ into the hemispheres, one of which
contains $i$ vertices from one component of $\widetilde{\G}$ and
$m-i$ vertices from the other component.
This gives relations $a^ib^{m-i}=a^{m-i}b^i=1$.

\enddocument